\documentclass{article}
\usepackage{graphicx}
\usepackage{psfrag}
\usepackage{amssymb}
\usepackage{amsmath}
\usepackage{amsthm}
\usepackage{amscd}

\newtheorem{thm}{Theorem}[section]
\newtheorem{prop}[thm]{Proposition}

\newtheorem{lem}[thm]{Lemma}

\newtheorem*{LS}{L$\acute{\rm e}$vy-Steinitz Theorem}
\newtheorem*{CT}{Coriolis Test}
\newtheorem*{rh}{Riemann Hypothesis}

\theoremstyle{definition}
\newtheorem{defn}[thm]{Definition}
\newtheorem{exmp}[thm]{Example}

\title{Ordinality and Riemann Hypothesis II}
\author{Young Deuk Kim
\\SNU College\\
Seoul National University\\
Seoul 08826, Korea
\\(yk274@snu.ac.kr)}
\date{\today}

\begin{document}
\maketitle

\begin{abstract}
For $\frac{1}{2}<x<1$, $y>0$, and $n\in\mathbb{N}$, let $\displaystyle\theta_n(x+iy)=\sum_{i=1}^n\frac{{\mbox{sgn}}\, q_i}{q_i^{x+iy}}$,
where $Q=\{q_1,q_2,q_3,\cdots\}$ is the set of finite products of distinct odd primes, and
${\mbox{sgn}}\, q=(-1)^k$ if $q$ is the product of $k$ distinct primes.
In this paper, we prove that there exists an ordering of $Q$ such that the sequence $\theta_n(x+iy)$ has a convergent subsequence.
As an application, we study the Riemann hypothesis.\\
\vspace{.3cm}
\noindent
2020 Mathematics Subject Classification ; 11M26.
\end{abstract}

\section{Introduction}

Let $\mathbb{N}$ be the set of natural numbers and $P$ be the set of odd primes.
\begin{defn}
For an ordering of $P=\{p_1,p_2,p_3,\cdots\}$ and $m\in\mathbb{N}$, let
\begin{equation*}
P_m=\{p_1,p_2,\cdots,p_m \}.
\end{equation*}
\end{defn}

\begin{defn}\label{defn:um}
Let $Q$ be the set of finite products of distinct odd primes.
\begin{equation*}
Q=\{ p_1p_2 \cdots p_k \mid k\in\mathbb{N} \mbox{ and } p_1, p_2,\cdots, p_k \mbox{ are distinct primes in }P\}
\end{equation*}
and, for each $m\in\mathbb{N}$, let
\begin{equation*}
U_m=\{ p_1p_2 \cdots p_k \mid k\in\mathbb{N} \mbox{ and }p_1, p_2,\cdots, p_k \mbox{ are distinct primes in }P_m\}.
\end{equation*}
Note that $U_m$ depends on the choice of ordering of $P$, and $U_m\subset U_{m+1}$.
\end{defn}

\begin{lem}\label{lem:um}
The number of elements in $U_m$ is $2^m-1$.
\end{lem}
\begin{proof}
Since
$$U_m=\{p_1,\cdots, p_m,\ p_1p_2,\cdots,p_{m-1}p_m,\ p_1p_2p_3,\cdots\cdots, p_1p_2\cdots p_m\},$$
the number of elements in $U_m$ is
$$\binom{m}{1}+\binom{m}{2}+\cdots +\binom{m}{m}=2^m-1.$$
\end{proof}

\begin{defn}
Let
\begin{equation*}
Q_1=U_1 \mbox{ and }Q_m=U_m-U_{m-1} \mbox{ for each }m=2,3,4, \cdots.
\end{equation*}
\end{defn}

\noindent
Notice that
\begin{equation}\label{eq:important}
Q_m=\{p_m,p_mq\mid q\in U_{m-1}\},\qquad\quad  \bigcup_{i=1}^m Q_m=U_m
\end{equation}
and $Q_1,Q_2,Q_3,\cdots$ are mutually disjoint.
Observe also that the number of elements in $Q_m$ is
\begin{equation*}
(2^m-1)-(2^{m-1}-1)=2^{m-1}.
\end{equation*}

\begin{exmp}\label{exmp;Q}
In the increasing ordering of $P$, we have
$$p_1=3,\ p_2=5,\ p_3=7,\cdots.$$
Therefore
\begin{eqnarray*}
Q_1=\{3\},\ Q_2=\{5,\ 3\cdot5\},\ Q_3=\{7,\ 3\cdot 7,\ 5\cdot 7,\ 3\cdot 5\cdot 7 \}, \cdots.
\end{eqnarray*}
\end{exmp}

\begin{defn}\label{defn:ordering}
An ordering of $P$, together with the following two conditions (C1)-(C2), induce a unique ordering of $Q=\{q_1,q_2,q_3\cdots \}$.
\begin{enumerate}
\item[(C1)]
$i<j$ if  $q_i<q_j \mbox{ and } q_i,q_j\in Q_m \mbox{ for some }m$.
\item[(C2)]
$i<j$ if $q_i\in Q_m,\ q_j\in Q_n\mbox{ for some }m<n$
\end{enumerate}
\end{defn}

Note that any ordering of $P$ induces a unique ordering of $Q$ in this way.

\begin{exmp} Suppose that $P$ has the increasing ordering. In the induced ordering of $Q$, we have
\begin{eqnarray*}
q_1=3,\ q_2=5,\ q_3=15,\ q_4=7,\ q_5=21,\ q_6=35,\ q_7=105,\ q_8=11,\cdots.
\end{eqnarray*}
\end{exmp}

\begin{defn}
For each $q= p_1p_2 \cdots p_k\in Q$, let $\mbox{sgn\,}q=(-1)^k$,
where $p_1, p_2,\cdots, p_k$ are distinct odd primes.
\end{defn}

\begin{defn}\label{def:theta}
Suppose that an ordering is given on $Q=\{q_1,q_2,q_3,\cdots\}$.
For $\frac{1}{2}<x<1$, $y>0$ and $n\in\mathbb{N}$, let
$$\theta_n(x+iy)=\sum_{i=1}^n\frac{\mbox{sgn\,}q_i}{q_i^{x+iy}}$$
\end{defn}

In this paper, we prove

\begin{thm}\label{mainthm}
For each $\frac{1}{2}<x<1$ and $y>0$, there exists an ordering of $P$ such that, under the induced ordering of $Q$,
$\theta_n(x+iy)$ has a convergent subsequence.
\end{thm}

As an application, we study the Riemann hypothesis.

\section{Preliminary Theorems}

We need the following theorem in the proof of Theorem \ref{mainthm}.

\begin{thm}[\cite{kim}]\label{kim}
Suppose that $y>0$, $0\leq\alpha<2\pi$ and $0<K<1$. Let $P^+$ be the set of primes $p$ such that $\cos(y\ln p+\alpha)>K$ and $P^-$ the set of primes $p$ such that $\cos(y\ln p+\alpha)<-K$ .
Then we have
$$\sum_{p\in P^+}\frac{1}{p}=\infty\quad\mbox{and}\quad \sum_{p\in P^-}\frac{1}{p}=\infty.$$
\end{thm}

From the argument in the proof of the Riemann rearrangement theorem, we have
\begin{thm}[\cite{riemann}, \cite{rosenthal}]\label{rosenthal}
For a series $\sum_{i=1}^\infty a_i$ of real numbers, suppose that
$$\lim_{i\to\infty}a_i=0$$
and let
\begin{equation}\label{eq:pm}
a_i^+=\mbox{max\,}\{a_i,0\}\qquad\mbox{and}\qquad a_i^-=-\mbox{min\,}\{a_i,0\}.
\end{equation}
If
$$\sum_{i=1}^\infty a_i^+=\sum_{i=1}^\infty a_i^-=\infty$$
then there exists a rearrangement such that the series $\sum_{i=1}^\infty a_i$ is convergent.
\end{thm}

We need the L$\acute{e}$vy-Steinitz theorem which is a generalization of the Riemann rearrangement theorem and Theorem \ref{rosenthal}.
\begin{LS}[\cite{rosenthal}]\label{Rosental}
The set of all sums of rearrangements of a given series of vectors
$$\sum_{i=1}^\infty \textbf{v}_i$$
in $\mathbf{R}^n$ is either the empty set or a translate of a subspace i.e., a set of the form $\mathbf{v}+M$, where $\mathbf{v}$ is a vector
and $M$ is a subspace. If the following two conditions (a)-(b) are satisfied then it is nonempty i.e., it has convergent
rearrangements.
\begin{enumerate}
\item[(a)]
$\lim_{i\to\infty}\mathbf{v}_i=\mathbf{0}$
\item[(b)]
For all vector $\mathbf{w}$ in $\mathbb{R}^n$,
$$\sum_{i=1}^\infty (\mathbf{v}_i,\mathbf{w})^+\quad\mbox{and}\quad \sum_{i=1}^\infty (\mathbf{v}_i,\mathbf{w})^-$$
are either both finite or both infinite, where we use the notations in eq. (\ref{eq:pm}) and
$(\mathbf{v}_i,\mathbf{w})$ is the Euclidean inner product of $\mathbf{v}_i$ and $\mathbf{w}$.
\end{enumerate}
\end{LS}

The Coriolis test will be used in the proof of Theorem \ref{mainthm}.
\begin{CT}[\cite{wermuth}]\label{wermuth}
If $z_i$ is a sequence of complex numbers such that
$$\sum_{i=1}^\infty z_i\quad \mbox{and}\quad \sum_{i=1}^\infty |z_i|^2$$
are convergent, then
$$\prod_{i=1}^\infty (1+z_i)$$
converges.
\end{CT}

\section{Proof of Theorem \ref{mainthm}}

\begin{defn}
Suppose that $P$ has the increasing ordering.
For $0<x\leq 1$ and $y\in\mathbb{R}$, let
\begin{eqnarray*}
\rho(x+iy)&=&\frac{1}{2^{x+iy}}+\sum_{i=1}^\infty\frac{1}{p_i^{x+iy}}\\
&=&\frac{\cos(y\ln 2)-i\sin(y\ln 2)}{2^x}+\sum_{i=1}^\infty\frac{\cos(y\ln p_i)-i\sin(y\ln p_i)}{p_i^x}
\end{eqnarray*}
\end{defn}

\begin{thm}\label{thm:rho}
Let $0<x\leq 1$ and $y\in\mathbb{R}$.
$\rho(x+iy)$ has a convergent rearrangement, and therefore
\begin{equation}\label{eqr}
\sum_{i=1}^\infty\frac{1}{p_i^{x+iy}}
\end{equation}
has a convergent rearrangement, too. In other words, $P$ has an ordering such that
eq. (\ref{eqr}) is convergent.
\end{thm}

\begin{proof}
Let
$$\mathbf{v}_1=\left(\frac{\cos(y\ln 2)}{2^x},\ -\frac{\sin(y\ln 2)}{2^x}\right)$$
and, for $i\in\mathbb{N}$, let
$$\mathbf{v}_{i+1}=\left(\frac{\cos(y\ln p_i)}{p_i^x},\ -\frac{\sin(y\ln p_i)}{p_i^x}\right).$$
Since $P$ has the increasing ordering, we have
\begin{equation}\label{eq:limit}
\lim_{i\to\infty}\mathbf{v}_i=\mathbf{0}.
\end{equation}

\indent
Let $$\mathbf{w}=r(\cos\alpha,\ \sin\alpha)$$ be a vector in $\mathbb{R}^2$, where $r\geq 0$ and $0\leq\alpha<2\pi$.
If $r=0$ then $(\mathbf{v}_i,\mathbf{w})=0$ for all $i\in\mathbb{N}$ and therefore
\begin{equation}\label{eq0}
\sum_{i=1}^\infty (\mathbf{v}_i,\mathbf{w})^+=\sum_{i=1}^\infty (\mathbf{v}_i,\mathbf{w})^-=0.
\end{equation}

\noindent
Suppose that $r>0$. We have
\begin{eqnarray*}
\mathbf{v}_1\cdot\mathbf{w}&=&\frac{r\cos(y\ln 2)\cos\alpha-r\sin(y\ln 2)\sin\alpha}{2^x}\\
&=&\frac{r\cos(y\ln 2+\alpha)}{2^x}
\end{eqnarray*}
and
\begin{eqnarray*}
\mathbf{v}_{i+1}\cdot\mathbf{w}&=&\frac{r\cos(y\ln p_i)\cos\alpha-r\sin(y\ln p_i)\sin\alpha}{p_i^x}\\
&=&\frac{r\cos(y\ln p_i+\alpha)}{p_i^x}
\end{eqnarray*}

\noindent
Let $P^+$ be the set of primes $p$ such that $\cos(y\ln p+\alpha)>\frac{1}{2}$ and $P^-$ the set of primes $p$
such that $\cos(y\ln p+\alpha)<-\frac{1}{2}$. From Theorem \ref{kim}, we have

\begin{eqnarray*}
\sum_{i=1}^\infty (\mathbf{v}_i,\mathbf{w})^+&\geq& \sum_{p\in P^+}\frac{r\cos(y\ln p+\alpha)}{p^x}\\
&& \geq\frac{r}{2}\sum_{p\in P^+}\frac{1}{p^x}\geq\frac{r}{2}\sum_{p\in P^+}\frac{1}{p}=\infty
\end{eqnarray*}
and
\begin{eqnarray*}
\sum_{i=1}^\infty (\mathbf{v}_i,\mathbf{w})^-&\geq& -\sum_{p\in P^-}\frac{r\cos(y\ln p+\alpha)}{p^x}\\
&& \geq\frac{r}{2}\sum_{p\in P^-}\frac{1}{p^x}\geq\frac{r}{2}\sum_{p\in P^-}\frac{1}{p}=\infty.
\end{eqnarray*}
Therefore
\begin{equation}\label{eq:infty}
\sum_{i=1}^\infty (\mathbf{v}_i,\mathbf{w})^+=\sum_{i=1}^\infty (\mathbf{v}_i,\mathbf{w})^-=\infty.
\end{equation}
\indent
From eq. (\ref{eq:limit}), (\ref{eq0}), (\ref{eq:infty}) and L$\acute{\rm e}$vy-Steinitz theorem, we know that the series of vectors in $\mathbb{R}^2$
$$\sum_{i=1}^\infty \textbf{v}_i$$
has a convergent rearrangement, and therefore $\rho(x+iy)$ has a convergent rearrangement.
\end{proof}

\begin{lem}\label{lem:PS}
Let $z=x+iy$, where $x,y\in\mathbb{R}$. For all $m\in\mathbb{N}$, we have
$$\prod_{i=1}^m\left(1-\frac{1}{p_i^z}\right)-1=\sum_{q\in Q_1}\frac{{\mbox{\rm{sgn}}}\, q}{q^z}+\sum_{q\in Q_2}\frac{{\mbox{\rm{sgn}}}\, q}{q^z}+\cdots+
\sum_{q\in Q_m}\frac{{\mbox{\rm{sgn}}}\, q}{q^z}.$$
\end{lem}
\begin{proof}
We use induction on $m$. If $m=1$, it is clear. Suppose that it is true for $m=k-1$.
We will show that it is true for $m=k$.
From eq. (\ref{eq:important}), we have
\begin{eqnarray*}
\prod_{i=1}^k\left(1-\frac{1}{p_i^z}\right)&=&\left(\prod_{i=1}^{k-1}\left(1-\frac{1}{p_i^z}\right)\right)\left(1-\frac{1}{p_k^z}\right)\\
&=&\left(1+\sum_{q\in Q_1}\frac{{\mbox{sgn}}\, q}{q^z}+\cdots+\sum_{q\in Q_{k-1}}\frac{{\mbox{sgn}}\, q}{q^z}\right)\left(1-\frac{1}{p_k^z}\right)\\
&=&\left(1+\sum_{q\in Q_1}\frac{{\mbox{sgn}}\, q}{q^z}+\cdots+\sum_{q\in Q_{k-1}}\frac{{\mbox{sgn}}\, q}{q^z}\right)\\
&&\qquad\qquad\quad  -\frac{1}{p_k^z}\left(1+\sum_{q\in Q_1}\frac{{\mbox{sgn}}\, q}{q^z}+\cdots+\sum_{q\in Q_{k-1}}\frac{{\mbox{sgn}}\, q}{q^z}\right)\\
&=&\left(1+\sum_{q\in Q_1}\frac{{\mbox{sgn}}\, q}{q^z}+\cdots+\sum_{q\in Q_{k-1}}\frac{{\mbox{sgn}}\, q}{q^z}\right)-\frac{1}{p_k^z}\left(1+\sum_{q\in U_{k-1}}\frac{{\mbox{sgn}}\, q}{q^z}\right)\\
&=&1+\sum_{q\in Q_1}\frac{{\mbox{sgn}}\, q}{q^z}+\cdots+\sum_{q\in Q_{k-1}}\frac{{\mbox{sgn}}\, q}{q^z}+
\sum_{q\in Q_k}\frac{{\mbox{sgn}}\, q}{q^z}
\end{eqnarray*}
\end{proof}

\indent
We are now ready to prove Theorem \ref{mainthm}.\\

\noindent
{\textbf{\large{Proof of Theorem \ref{mainthm}}}}\\

Recall that $\frac{1}{2}<x<1$ and $y>0$.
By Theorem \ref{thm:rho}, we can choose an ordering of $P$ such that
$$\sum_{i=1}^\infty\frac{1}{p_i^{x+iy}}$$ is convergent.
From now on, we assume that $P$ has the chosen ordering, and $Q$ has the induced ordering.

Since $\frac{1}{2}<x<1$,
$$\sum_{i=1}^\infty \left| \frac{1}{p_i^{x+iy}}\right|^2 =\sum_{i=1}^\infty\frac{1}{p_i^{2x}}$$
is convergent. Therefore, by the Coriolis test,
$$\prod_{i=1}^\infty \left(1-\frac{1}{p_i^{x+iy}}\right)$$
is convergent. By Lemma \ref{lem:PS}, Lemma \ref{lem:um} and eq. (\ref{eq:important}),  we have
\begin{eqnarray*}
\prod_{i=1}^m\left(1-\frac{1}{p_i^{x+iy}}\right)-1&=&\sum_{q\in Q_1}\frac{{\mbox{sgn}}\, q}{q^{x+iy}}+\sum_{q\in Q_2}\frac{{\mbox{sgn}}\, q}{q^z}+\cdots+\sum_{q\in Q_m}\frac{{\mbox{sgn}}\, q}{q^{x+iy}}\\
&=&\sum_{q\in U_m}\frac{{\mbox{sgn}}\, q}{q^{x+iy}}\\
&=&\sum_{i=1}^{2^m-1}\frac{{\mbox{sgn}}\, q_i}{q_i^{x+iy}}.
\end{eqnarray*}
Therefore
\begin{equation}\label{eq:consub}
\theta_{2^m-1}(x+iy)=\sum_{q\in U_m}\frac{{\mbox{sgn}}\, q}{q^{x+iy}}
\end{equation}
is a convergent subsequence of $\theta_n(x+iy)$.
\qed

\section{Application to the Riemann Hypothesis}

The zeta function was introduced by Euler in 1737 for real variable $s>1$.
In 1859, Riemann \cite{riemann2} extended the function to a complex meromorphic function $\zeta(z)$
with only simple pole at $z=1$.

\begin{rh}[\cite{bombieri}, \cite{sarnak}]
The zeros of $\zeta(z)$ in the critical strip  $0<{Re\,}z<1$ lie on the critical line $\mbox{Re\,} z=\frac{1}{2}$.
\end{rh}

\indent
Suppose that $x$ and $y$ are real numbers with $0<x<1$. It is known that if $x+yi$ is a zero of the zeta function, then so are $x-yi$, $(1-x)+yi$, and $(1-x)-yi$.

\indent
Riemann himself showed that if $0<x<1$, $0\leq y\leq 25.02$ and $x+yi$ is a zero of the zeta function, then $x=\frac{1}{2}$. Therefore the Riemann hypothesis is true up to height 25.02.
In 1986, van de Lune, te Riele and Winter \cite{LRW} showed that the Riemann hypothesis is true up to height 545,439,823,215. Furthermore in 2021 Dave Platt and Tim Trudgian \cite{PT} proved that the Riemann hypothesis is true up to height $3\cdot 10^{12}$.

\indent
Therefore, to prove the Riemann hypothesis, it is enough to show that if $\frac{1}{2}<x<1$ and $y>0$, then $x+yi$ is not a zero of the zeta function.

The eta function
$$\eta(z)=\sum_{k=1}^\infty\frac{(-1)^{k-1}}{k^z}=1-\frac{1}{2^z}+\frac{1}{3^z}-\cdots$$
is convergent on $\operatorname{Re}(z)>0$ and is useful in the study of the Riemann hypothesis.

\begin{thm}[\cite{broughan}]\label{thm:ext}
In $0<\operatorname{Re}(z)<1$, any zero of $\zeta(z)$ is a zero of $\eta(z)$.
\end{thm}

Therefore, to prove the Riemann hypothesis, it is enough to show that if $\frac{1}{2}<x<1$ and $y>0$, then $x+yi$ is not a zero of the eta function.

\indent
We will use the following notation.
\begin{defn}\label{defn:phi}
Suppose that $0<x<1$ and $y\in\mathbb{R}$ are given. For each $k\in\mathbb{N}$, let
$$\varphi(k)=\frac{(-1)^{k-1}}{k^{x+iy}},$$
where we assume that $(-1)^0=1$ for the sake of simplicity.
\end{defn}

\indent
The following theorem is crucial in the proof of the Riemann hypothesis.
We include the proof for completeness.
\begin{thm}[\cite{kim2}]\label{thm:2n}
Suppose that $0<x<1$ and $y\in\mathbb{R}$. Then
\begin{equation*}
\sum_{\ell=0}^\infty\varphi(2^\ell)
\end{equation*}
converges to a nonzero complex number.
\end{thm}

\begin{proof}
Since $0<x<1$, we have
\begin{eqnarray*}
&&\sum_{\ell=0}^\infty\varphi(2^\ell)=1-\sum_{\ell=1}^\infty\frac{1}{(2^\ell)^{x+iy}}=1-\sum_{\ell=1}^\infty\frac{e^{-i\ell y\ln 2}}{2^{\ell x}}\\
&&\qquad\quad =1-\sum_{\ell=1}^\infty \left(\frac{e^{-iy\ln 2}}{2^x}\right)^\ell = 1-\frac{e^{-iy\ln 2}}{2^x-e^{-iy\ln 2}}=\frac{2^x-2e^{-iy\ln 2}}{2^x-e^{-iy\ln 2}}\neq 0.
\end{eqnarray*}
\end{proof}

\begin{defn}
Suppose that an ordering is given on $Q=\{q_1,q_2,q_3,\cdots\}$. For $k,i\in\mathbb{N}$, let
\[\delta(k,i)=\left\{
\begin{array}{cl}
1 & \mbox{ if } k \mbox{ is a multiple of } q_i \\
0& \mbox{ otherwise}
\end{array}\right.
\]
\end{defn}

\begin{defn}
Let $\mathbb{N}_0=\mathbb{N}\cup\{ 0\}$.
\end{defn}

\begin{defn}\label{defn:fk}
For $k\in\mathbb{N}$, let
\begin{eqnarray*}
f(k)&=&\left\{
\begin{array}{cl}
0& \mbox{if } k=2^\ell \mbox{ for some }\ell\in\mathbb{N}_0\\
-1 & \mbox{otherwise}
\end{array}\right.
\end{eqnarray*}
\end{defn}

The author's research on the Riemann hypothesis is motivated by the following theorem.
We include its proof for the sake of completeness.

\begin{thm}[\cite{kim2}]\label{thm:fk}
Let $k\in\mathbb{N}$. For any ordering of $Q=\{q_1,q_2,q_3,\cdots\}$, we have
$$\sum_{i=1}^\infty(\mbox{\rm{sgn}}\;q_i)\delta(k,i)=f(k)$$
\end{thm}

\begin{proof}
If $k=2^\ell$ for some $\ell\in\mathbb{N}_0$, then $k$ is not a multiple of any element in $Q$. Therefore
$\delta(k,i)=0$ for all $i\in\mathbb{N}$ and hence
$$\sum_{i=1}^\infty(\mbox{\rm{sgn}}\;q_i)\delta(k,i)=0=f(k).$$

\indent
Suppose that $k\neq 2^\ell$ for any $\ell\in\mathbb{N}_0$. Let
$$k=2^mp_1^{m_1}p_2^{m_2}\cdots p_n^{m_n}$$
be the prime factorization of $k$, where $p_1,p_2,\cdots, p_n$ are distinct odd primes.
Notice that $n,m_1, m_2,\cdots, m_n\in\mathbb{N}$ and $m\in\mathbb{N}_0$.
We have
\begin{eqnarray*}
&&\{q_i\in Q\mid \delta(k,i)=1\}\\
&=&\{p_1,\cdots, p_n,\ p_1p_2,\cdots,p_{n-1}p_n,\ p_1p_2p_3,\cdots\cdots, p_1p_2\cdots p_n\}.
\end{eqnarray*}
Therefore
\begin{eqnarray*}
\sum_{i=1}^\infty(\mbox{\rm{sgn}}\;q_i)\delta(k,i)&=&-\binom{n}{1}+\binom{n}{2}-\cdots +(-1)^n\binom{n}{n}\\
&=&(1-1)^n-1=f(k).
\end{eqnarray*}
\end{proof}

\subsection{Proof of the Riemann hypothesis}

\noindent
In this section, we assume that $z=x+iy$ is a zero of $\eta(z)$, where $\frac{1}{2}<x<1$ and $y>0$.
This leads to a contradiction, thereby proving the Riemann hypothesis.

\indent
By Theorem \ref{thm:rho}, we can choose an ordering of
$P=\{p_1,p_2,p_3,\cdots\}$ such that
$$\sum_{i=1}^\infty\frac{1}{p_i^{x+iy}}$$ is convergent.
From now on, we assume that $$P=\{p_1,p_2,p_3,\cdots\}$$ has the chosen ordering and
$$Q=\{q_1,q_2,q_3,\cdots\}$$ has the induced ordering.
Notice that
$$\lim_{m\to\infty}\sum_{i=1}^{2^m-1}\frac{\mbox{\rm{sgn}}\;q_i}{q_i^z}$$
is convergent. Notice also that
$$\sum_{h=1}^\infty\frac{(-1)^{h-1}}{h^z}=0.$$

\begin{defn}
For $m,n\in\mathbb{N}$, let
$$\Phi(m,n)=\left(\sum_{i=1}^{2^m-1}\frac{\mbox{\rm{sgn}}\;q_i}{q_i^z}\right)\left(\sum_{h=1}^n\frac{(-1)^{h-1}}{h^z}\right)=
\left(\sum_{q_i\in U_m}\frac{\mbox{\rm{sgn}}\;q_i}{q_i^z}\right)\left(\sum_{h=1}^n\frac{(-1)^{h-1}}{h^z}\right)$$
\end{defn}

\indent
Since any $q_i\in Q$ is an odd number, we have
  $$(-1)^{h-1}=(-1)^{hq_i-1}$$
for all $i$. Let $k=hq_i$. We have
\begin{eqnarray}
\Phi(m,n)&=&\sum_{h=1}^n\sum_{q_i\in U_m}\frac{(\mbox{\rm{sgn}}\;q_i)(-1)^{hq_i-1}}{(hq_i)^z} \nonumber\\
&=&\sum_{k=1}^{np_1p_2\cdots p_m} \left(\sum_{q_i\in U_m}^{q_i\geq k/n}(\mbox{\rm{sgn}}\;q_i)\delta(k,i)\right)\frac{(-1)^{k-1}}{k^z} \label{eq:Phi}
\end{eqnarray}

\begin{defn}
Let $K\in\mathbb{N}$. Let $m(K)$ be the smallest natural number $m$ such that
$U_m$ contains all odd primes less than or equal to K.
\end{defn}

\begin{defn}\label{defn:Psi}
Let
\begin{eqnarray*}
\Psi(K)=\Phi(m(K), K)=\left(\sum_{q_i\in U_{m(K)}}\frac{\mbox{\rm{sgn}}\;q_i}{q_i^z}\right)\left(\sum_{h=1}^K\frac{(-1)^{h-1}}{h^z}\right).
\end{eqnarray*}
\end{defn}

\noindent
Notice that
\begin{equation}\label{eq:psi}
\lim_{K\to\infty}\Psi(K)=0.
\end{equation}

\indent
We have
\begin{eqnarray}
\Psi(K)&=&\sum_{h=1}^K\sum_{q_i\in U_{m(K)}}\frac{(\mbox{\rm{sgn}}\;q_i)(-1)^{hq_i-1}}{(hq_i)^z} \nonumber\\
&=&\sum_{k=1}^{Kp_1p_2\cdots p_{m(K)}} \left(\sum_{q_i\in U_{m(K)}}^{q_i\geq k/K}(\mbox{\rm{sgn}}\;q_i)\delta(k,i)\right)\frac{(-1)^{k-1}}{k^z} \nonumber\\
&=&\sum_{k=1}^{Kp_1p_2\cdots p_{m(K)}} c(K,k)\varphi(k),\label{eq:Psi}
\end{eqnarray}
where
\begin{eqnarray}\label{eq:cmk1}
c(K,k)=\sum_{q_i\in U_{m(K)}}^{q_i\geq k/K}(\mbox{\rm{sgn}}\;q_i)\delta(k,i).
\end{eqnarray}

For all $1\leq k\leq K$, from Theorem \ref{thm:fk}, we have
\begin{equation}\label{eq:cKk}
c(K,k)=\sum_{q_i\in U_{m(K)}}(\mbox{\rm{sgn}}\;q_i)\delta(k,i)=\sum_{i=1}^\infty(\mbox{\rm{sgn}}\;q_i)\delta(k,i)=f(k).
\end{equation}

Therefore, from eq. (\ref{eq:Psi}) and eq. (\ref{eq:cKk}), we have
\begin{eqnarray}
\Psi(K)&=&\sum_{k=1}^{Kp_1p_2\cdots p_{m(K)}} c(K,k)\varphi(k) \nonumber\\
&=&\sum_{k=1}^K c(K,k)\varphi(k)+\sum_{k=K+1}^{Kp_1p_2\cdots p_{m(K)}} c(K,k)\varphi(k)\nonumber\\
&=&\sum_{k=1}^K f(k)\varphi(k)+\sum_{k=K+1}^{Kp_1p_2\cdots p_{m(K)}} c(K,k)\varphi(k).\label{eq:Psi2}
\end{eqnarray}

\begin{defn}\label{defn:Omega}
Let
$$\Omega(K)=\sum_{k=K+1}^{Kp_1p_2\cdots p_{m(K)}} c(K,k)\varphi(k).$$
\end{defn}

\begin{prop}\label{prop:difficult}
We have
$$\lim_{K\to\infty}\Omega(K)=0$$
and therefore
$$\lim_{K\to\infty}\Psi(K)=\sum_{k=1}^\infty f(k)\varphi(k).$$
\end{prop}

\begin{proof}
The proof of this proposition in the previous versions is wrong.
I will try to find a correct proof of this proposition.
\end{proof}

\indent
Suppose that Proposition \ref{prop:difficult} is true. From eq. (\ref{eq:psi}), we have
\begin{equation}\label{eq:victory}
\sum_{k=1}^\infty f(k)\varphi(k)=\lim_{K\to\infty}\Psi(K)=0.
\end{equation}

\noindent
Since $z=x+iy$ is a zero of the $\eta(z)$, we have
$$\sum_{k=1}^\infty \varphi(k)=0.$$
Thus, from eq. (\ref{eq:victory}), we have
$$\sum_{\ell=0}^\infty \varphi(2^\ell)=\sum_{k=1}^\infty \varphi(k)+\sum_{k=1}^\infty f(k)\varphi(k)=0.$$
This contradicts Theorem \ref{thm:2n}.
Therefore, the Riemann hypothesis is true if Proposition \ref{prop:difficult} is true.



\begin{thebibliography}{00}

\bibitem {bombieri} E. Bombieri, \textit{Problems of the Millennium: The Riemann Hypothesis}.

\bibitem {broughan} K. Broughan, \textit{Equivalents of the Riemann Hypothesis Volume One: Arithmetic Equivalents},
Encyclopedia of Mathematics and Its Applications 164, Cambridge University Press, 2017.

\bibitem {kim} Y. D. Kim, \textit{On the sum of reciprocals of primes}, preprint, arXiv:2403.04768.

\bibitem {kim2} Y. D. Kim, \textit{Ordinality and Riemann Hypothesis I}, preprint, arXiv:2311.00003.

\bibitem {LRW} J. van de Lune, H. J. J. te Riele and D. T. Winter, \textit{On the zeros of
the Riemann zeta function in the critical strip IV}, Math. Comp. {\bf 46}(1986), 667-681.

\bibitem {riemann} B. Riemann, \textit{\"{U}ber die Darstellbarkeit einer Function durch eine trigonometrische Reihe}, Abh. kgl. Ges. Wiss. Göttingen 13, 87–132 (1867) = Gesammelte mathematische Werke (Leipzig 1876), 213-253.

\bibitem {riemann2} B. Riemann, \textit{\"{U}ber die Anzahl der Primzahlen unter einer gegebenen Gr\"{o}sse},
Monatsberichte der Berliner Akademie (1859), 671-680.

\bibitem {rosenthal} P. Rosenthal, \textit{A remarkable Theorem of Levy and Steinitz}, The American Mathematical Monthly {\bf 94}(4) (1987), 342-351.

\bibitem {PT} D. Platt, T. Trudgian, \textit{The Riemann hypothesis is true up to $3\cdot 10^{12}$}, Bulletin of the London Mathematical Society {\bf 53}(2021), 792-797.

\bibitem {sarnak} P. Sarnak, \textit{Problems of the Millennium: The Riemann Hypothesis}, CLAY 2004.

\bibitem {wermuth} E. Wermuth, \textit{Some Elementary Properties of Infinite Products}, The American Mathematical Monthly {\bf 99}(6) (1992), 530-537.

\end{thebibliography}
\end{document}